\newcommand{\C}{{\mathbf C}}
\newcommand{\Q}{{\mathbf Q}}
\newcommand{\Qbar}{\overline{\Q}}
\newcommand{\Z}{{\mathbf Z}}
\newcommand{\N}{{\mathbf N}}
\newcommand{\p}{{\mathbf P}}
\newcommand{\F}{{\mathbf F}}
\newcommand{\R}{{\mathbf R}}
\newcommand{\hs}{{\mathbf H}}
\newcommand{\hh}{{\cal H}}
\newcommand{\oo}{{\cal O}}
\newcommand{\pp}{{\cal P}}
\newcommand{\half}{\frac 12}
\newcommand{\disc}{\mbox{\rm \,disc\,}}
\newcommand{\SL}{\mbox{\rm \,SL\,}}

\documentclass{llncs}
\begin{document}
%
%
\title{
Explicit models of genus 2 curves with split CM
}
\titlerunning{Explicit Models}  
%
\author{
 Fernando Rodriguez-Villegas
}
\authorrunning{
 Fernando Rodriguez-Villegas
}   
%
\tocauthor{
 Fernando Rodriguez-Villegas
}
\institute{University of Texas at Austin, TX 78712, USA,\\
\email{
villegas@math.utexas.edu},\\
 WWW home page:
\texttt{
http://www.ma.utexas.edu/users/villegas
}
}
\thanks{
Much of this work was done while I was a guest at the Max Planck
Institut f\"ur Mathematik in Bonn in 1995. I take the opportunity to
thank everybody at the Institut for their hospitality.
I would also like to thank the NSF, TARP, and the Alfred P.~Sloan
Foundation for their generous support.
}

\maketitle              

\begin{abstract}
We outline a general algorithm for computing  an explicit model over a
number field of any curve of genus 2 whose (unpolarized) Jacobian is
isomorphic to the product of two elliptic curves with CM by the same
order in an imaginary quadratic field. We give the details and some
examples  for the case where the order has prime discriminant and
class number one.
\end{abstract}
\section{
Motivation
}

  Let $E_1, E_2$ be two elliptic curves defined over
$\Qbar$ with complex multiplication by an order $\oo$ of an imaginary
quadratic field $K$. We are interested in finding explicit models for
curves $C$ defined over $\Qbar$ whose (unpolarized) Jacobian is
isomorphic to $E_1 \times E_2$.  In this paper we propose a
general algorithm for this purpose and give details only for
the following special case where we have have carried them out:
$E_1=E_2$, $\oo$ is the ring of integers of $K=\Q(\sqrt{-N})$, $N\equiv
3 \bmod 4$ prime and $\oo$ has class number one.  Our special case
consists of finitely many curves, up to isomorphism; the algorithm
produces models over $K$ for them.

It is a general fact due to Narasimhan and Nori [NN] that there are
only finitely many principal polarizations on a given abelian variety
up to isomorphism.  Hence, for a fixed $\oo$ there are only finitely
many isomorphism classes of the curves we want; their number was
calculated by Hayashida and Nishi [HN].

For a similar question in the case of abelian surfaces with complex
multiplication by a quartic field see [vW].

 Our interest in this problem arose in connection with a
generalization to genus 2 of the {\it singular moduli} formulae of
Gross and Zagier [GZ] for the norm of the difference of $j$-values of
CM elliptic curves. (This generalization will be the subject of a
separate publication.) As an illustration, consider the genus 2 curve
$C$ determined by
$$
y^2=f(x)=6^{-3}h(x)h^\iota(x),
$$
where
\begin{eqnarray*}
h(x)=&(7144\sqrt{-163}-151790)x^3+(129789\sqrt{-163}+1752597)x^2\\
 &+(-47481\sqrt{-163}+510153)x +(-1596\sqrt{-163}-37250)\ ,
\end{eqnarray*}
and
$$
h^\iota(x)=\overline{x^3h(-1/x)}
$$
(bar denoting complex conjugation of the coefficients).
The unpolarized Jacobian of $C$ is isomorphic over $\Qbar$ to the
product of two 
elliptic curves with CM by the ring of integers of $K=\Q(\sqrt{-163})$. 
Let
$$
D=2^{-12}
\disc(f) 
=(2\cdot3^2\cdot5\cdot7\cdot11\cdot17\cdot19\cdot23)^{12}\ .
$$

Then we have
$$
\log D= -6\sum_{m \in
\Z^3}\sum_{d|(163-Q(m))/4}\left(\frac{-163}{d}\right)\log d\ ,
$$
where
$$
Q(m)= m^t  \left(
\begin{array}{ccc}
24&4&6\\
4&55&1\\
6&1&83
\end{array}
\right) m\ 
$$
is a certain positive definite  ternary quadratic form of
level $163$ associated to $C$ 
and the (finite) sum is over $m\in \Z^3$ such that $(163-Q(m))/4$ is a
positive integer.  In particular every rational prime $l$ dividing $D$
is  smaller than $163/4$ and inert in $K$.

The significance of the number $D$ is that $C$ has bad reduction only
at primes dividing $D$. Note that over $\Qbar$ the Jacobian of $C$ has
good reduction everywhere but $C$ does not; at primes diving $D$, $C$
reduces to two elliptic curves crossing at a point.

Another source of interest in the problem is the fact, which I learned
from K. Lauter, that the reduction of the curves $C$ provides genus 2
curves over certain finite fields with maximal number of rational
points (see  \S  5 for an example). In this regard, the more
interesting problem is the analogous one for curves of genus 3 for
which we hope to exhibit in the near future an algorithm similar to
the one sketched here.

\section{
Outline of the algorithm
}

 We start by giving an outline of the main steps of the
general algorithm and then  give  details for our special case in the
next sections.

\vskip .5cm

{\bf Step 1.}\hskip .3cm Find  period matrices  for the polarized Jacobians.
\vskip .2cm
{\bf Step 2.}\hskip .3cm Given a non-split period matrix obtained in step 1
compute a model for the corresponding curve.
\vskip .5cm

The first step is purely algebraic and only requires computations with
rational numbers; it involves the calculation of representatives for
ideal classes of certain orders in a quaternion algebra (see [HN] for
more details). What we need to do is describe explicitly the finitely
many principal polarizations on $E_1 \times E_2$ up to equivalence.

The second step relies on the following  explicit version of
Torelli's theorem for curves of 
genus 2 due to Bolza and Klein. Let 
$$
\hh_2=\{ Z \in \C^{2\times 2}\;\; |\quad Z^t=Z, \quad \mbox{\rm Im} (Z) \;\;
\mbox{\rm   positive definite}\}
$$ 
be the Siegel upper-half space of rank $2$. Let $Z \in \hh_2$ be a
period matrix of a  principally polarized abelian
surface which is not the product of two elliptic curves with the
product polarization. A theorem of Torelli guarantees that $Z$
arises from a  curve of genus 2, unique up to isomorphism.
Here is a way of recovering the curve from the period matrix $Z$.

 Let $f_Z(u_1,u_2) \in \C[u_1,u_2]$ be the
leading term in the Taylor expansion of
$$
\prod_{ (\mu,\nu) \mbox{\rm\  odd}} \theta_{\mu,\nu}(u,Z)\ ,
\qquad u=(u_1,u_2)
$$
about the origin, where for  $\mu, \nu \in \{0,1\}$ 
$$
\theta_{\mu,\nu}
(u,Z):=\sum_{m\in \Z^2+\half\mu}e^{\pi i  m^tZm}e^{2\pi i
  m^t(u+\half\nu)}\ , \quad u=(u_1,u_2) \in \C^2, \quad Z\in \hh_2
$$
is the  theta function with characteristics (see [Mu]). Then the
canonically polarized Jacobian of  the 
hyperelliptic curve over $\C$ determined by the equation 
$$
y^2=f_Z(x,1)
$$
corresponds to $Z$. (There are six  theta functions with odd
characteristics and hence
$f_Z$ is a sextic, i.e.  homogeneous of degree $6$.)

The difficulty in applying this formula is to know how to normalize
the sextic $f_Z$ properly to guarantee that its coefficients are
algebraic integers as well as finding similar expressions for its
Galois conjugates. In general, this would be accomplished by an
application of Shimura's general reciprocity law.  We would obtain
rapidly convergent series giving the minimal polynomials of these
coefficients. Since the coefficients of the minimal polynomials of the
coefficients of the sextic are in $\Z$, truncating the series would
then allow us to compute them {\em exactly}. We show how this works
for our special case in \S 4.

\section{
Principal polarizations
}

 From now on we assume that $\oo$ is the ring of integers of
$K=\Q(\sqrt{-N}) \subset \C$, $N\equiv 3 \bmod 4$ prime and the class
number of $\oo$ is $1$. Hence,  $E_1=E_2=E$, with $E$ isomorphic
to $\C/\oo$ over $\C$.

The principal polarizations of $E \times E$ up to isomorphism
correspond to positive definite unimodular Hermitian forms of rank 2
over $\oo$ up to $GL_2(\oo)$-equivalence. In order to find a set of
representatives of these Hermitian forms we will exploit the happy
accident that since we assume $N\equiv 3 \bmod 4$ the quaternion
algebra $B=(\frac{-1,-N}{\Q})$ (up to isomorphism the unique
quaternion algebra over $\Q$ ramified only at $N$ and $\infty$)
contains $\Q(i)$.  This allows us to convert the question to that of
finding Hermitian forms over $\Z[i]$ of discriminant $-N$ up to
equivalence and this is quite simple. Here is how it works.

Consider in $B$ the order
$$
R= \Z +\Z i+\Z\half(1+j)+\Z i\half(1+j), \qquad i^2=-1,j^2=-N.
$$
 $R$ is a maximal order in $B$ with a
natural embedding of $\oo$ sending $\sqrt{-N}$ to $j$. 
The rank 2 unimodular Hermitian forms arising from polarizations
of $E\times E$ correspond to rank $1$ left $R$-modules. 

Since $R$ also has an embedding of $\Z[i]$ (sending $i$ to $i$) we may
associate to a left $R$-module a rank $2$ Hermitian form $\Phi$ over
$\Z[i]$. We can give $\Phi$ as a triple $(a,b,c)$ with $a,c\in
\Z_{>0}$ and $b \in \Z[i]$, where $\Phi(u,v)=2au\overline u + bu
\overline v +\overline b \overline u v+2cv \overline v$.  It is not
hard to see that this form has discriminant $b\overline b -4ac=-N$. 

 It will be more convenient to work with $SL_2$ rather than $GL_2$
equivalence and to avoid duplications we consider only forms
$\Phi=(a,b,c)$ with $b \equiv 1 \bmod 2$. The above discussion
establishes a 1-1 correspondence between principal polarizations on
$E\times E$, up to $SL_2(\oo)$-equivalence, and positive definite
binary Hermitian forms $\Phi=(a,b,c)$ over $\Z[i]$ of discriminant
$-N$, up to $SL_2(\Z[i])$-equivalence.

Let $\hs$ be the hyperbolic $3$-space
$$
\hs=\{w=(x,y,t) \in \R^3 \;|\; t>0\}\ ,
$$
which we will think as embedded in the Hamilton quaternion algebra $H$ 
by $(x,y,t) \mapsto x+iy+jt$ (here $i,j$ are the usual basis of $H$
with $i^2=j^2=-1$ and $ij=-ji$).  If $\left(\begin{array}{cc}
  a&b\\c&d\end{array}\right) \in \SL_2(\C)$ then it is not hard 
to check that 
$$
w \mapsto (aw+b)(cw+d)^{-1}
$$
sends $\hs$ to $\hs$ defining an action of $SL_2(\C)$ on $\hs$ and in 
particular,  an action of $SL_2(\Z[i])$. This last action has a very
simple fundamental domain, whose closure is given by $w=(x,y,t)  \in
\hs$  with $x^2+y^2+t^2\geq 1, x\leq 1/2, y\leq 1/2, 0\leq x+y$. 

We can associate to a form $\Phi$ the point $w=(b+\sqrt{N}j)/2a \in
\hs$ and we call $\Phi$ {\it reduced}  if $w$ lies in the
fundamental domain.  The action of $\SL_2(\Z[i])$ on $\hs$ mimics that
on Hermitian forms. Every form is $\SL_2(\Z[i])$-equivalent to a
unique reduced form.

 The situation is in fact very analogous to that 
of positive definite binary quadratic forms over $\Z$ and, as in that
case, it is easy to write an algorithm that lists all reduced forms
$\Phi$ of a given discriminant (we do not really need $N$ to be prime
or class number $1$ for this).  Here is a brief sketch (all of this is
classical going back to Hermite [He, I, p. 251]).

{\tt

\noindent
Input: $N \equiv 3 \bmod 4$

\bigskip

\noindent
For $0\leq r,s \leq \sqrt{N/2}$, $r$ odd, $s$ even

 \quad Set $m:=(r^2+s^2+N)/4$

 \quad For  $a | m$, $\max(r,s) \leq a\leq \sqrt{m}$

  \quad \quad Add $(a,r+is,m/a)$ to List
 
 \quad \quad  Add $(a,-r+si,m/a)$ to List unless 

\quad \quad  $a=m/a$ or  $r=a$ or $s=a$ or  $s=0$

\bigskip

\noindent
Output: List
}

\bigskip
\bigskip
\noindent
As an example, we give in table 1 the list of reduced forms of discriminant
$-163$. 
\begin{table}
\begin{center}
\caption{Reduced Hermitian forms over $\Z[i]$ of discriminant $-163$}
\begin{tabular}{|c|}
\hline
$(1, 1, 41)$\\
$(2, 1 + 2i, 21)$\\
$(3, \pm1 + 2i, 14)$\\
$(6, \pm1 + 2i, 7)$\\
$(4, \pm3 + 2i, 11)$\\
$(6, \pm5 + 2i, 8)$\\
$(5, \pm1 + 4i, 9)$\\
$(7, \pm5 + 6i, 8)$\\
\hline
\end{tabular}
\end{center}
\end{table}

In general, the number of $\Phi$'s is the {\it class number} $n$ of
$B$ [Ei], which can be given in terms of $N$ as follows (a formula
valid for any prime $N \equiv 3 \bmod 4$)
$$
n=\left\{
\begin{array}{cc}
 \frac{1}{12}(N+5), & \quad \left(\frac{-3}N\right)=+1\\
\\
\frac1{12}(N+13), & \quad \left(\frac{-3}N\right)=-1.
\end{array}
\right.
$$

Finally, given a $\Phi=(a,b,c)$ as above, $b=r+si$,  the matrix 
$$
Z_\Phi:= \frac1{2a}
\left( \begin{array}{cc}
 r+\sqrt{-N} & s\\s&-r+\sqrt{-N}
\end{array} \right)
\in \hh_2 
$$
is a   period matrix  corresponding to the
associated principal polarization on $E\times E$.

\section{
Bolza-Klein sextics
}

 The product polarization on $E \times
E$ corresponds to the reduced form $\Phi=(1,1,(N+1)/4)$ in the
principal class; hence, forms $\Phi$ not in the principal class
correspond to curves. 

Given a form $\Phi=(a,b,c)$ not in the principal class we define the
associated normalized Bolza--Klein sextic $f_\Phi$ as follows.
$$
f_\Phi(u_1,u_2):=\frac{1}{a^6|\eta((1+\sqrt{-N})/2)|^{24}} \;
f_{Z_\Phi}(u_1,u_2),
$$
where $\eta$ is Dedekind's eta function.
It satisfies the following properties.
\begin{itemize}
\item
The $SL_2(\C)$ class of $f_\Phi$ depends only on the $SL_2(\Z[i])$
class of $\Phi$.
\item
$f_\Phi$ has coefficients in $K$ and $a^6f_\Phi$ has coefficients in
$\oo$. 
\item
The Igusa invariants [Ig] of $f_\Phi$ are in $\Z$ and depend only on the
$SL_2(\Z[i])$-equivalence class of $\Phi$.
\end{itemize}

The genus 2 curve
$$
C_\Phi: \qquad y^2=f_\Phi(x,1)
$$
is then defined over $K$ and, over the algebraic closure $\overline K$
of $K$ in $\C$, its Jacobian is isomorphic to $E\times E$. 

 Given a form $\Phi=(a,b,c)$ let $\Phi^\iota:=(a,-\overline
b,c)$. Suppose both $\Phi$ and $\Phi^\iota$ are reduced. Then $\Phi$
and $\Phi^\iota$ are not $SL_2(\Z[i])$-equivalent but they are
(always) $GL_2(\Z[i])$-equivalent. The corresponding curves $C_\Phi,
C_{\Phi^\iota}$ are hence isomorphic over $\overline K$; note that
they are also complex conjugates of each other. Otherwise,
curves $C_\Phi$ corresponding to different reduced forms are
non-isomorphic. The involution $\iota$ has a natural counterpart on
the left $R$-ideals in $B$ and it turns out that the number of orbits
of $\iota$ is what is classically known as the {\em type number} of
$B$ [Ei].   Hence, there are $t-1$ isomorphism classes of curves with
Jacobian isomorphic to $E\times E$, where $t$ is the
type number of the quaternion algebra $B$ [HN].

Here is a table with the values of $n$ and $t$ for the primes $N$ we are
considering. 
\begin{table}
\begin{center}
\caption{Type and class number of the quaternion algebra $B$}
\begin{tabular}{|r|r|r|}
\hline
N&n&t\\
\hline
3&1&1\\
7&1&1\\
11&2&2\\
19&2&2\\
43&4&3\\
67&6&4\\
163&14&8\\
\hline
\end{tabular}
\end{center}
\end{table}

Note that for $N=3$ or $7$ we only have the product polarization and hence
there is no curve $C$ with unpolarized Jacobian isomorphic to $E\times
E$ in that case.

Given a curve $C$ defined over $\Qbar$ its {\it field of moduli} is
the field $F\subset \Qbar$ characterized by the property: For every
$\tau \in Gal(\Qbar/\Q)$, $C^\tau$ is isomorphic to $C$ if and only if
$\tau$ is the identity on $F$. Clearly isomorphic curves have the same
field of moduli. Notice that $F$ is the smallest field over which a
curve isomorphic to $C$ {\it could} be defined, but it is not in
general a field over which it {\it can} be defined. In fact, for
example, Shimura showed that no generic hyperelliptic curve of even
genus has a model over its field of moduli [Sh Thm 3]. See [Me] for a
discussion of this issue for curves of genus 2.

For the curves $C_\Phi$ the field of moduli is $\Q$ (the field
generated by the Igusa invariants [Ig]), but, in fact, most are not
definable over $\Q$; only those forms $\Phi$ which are
$SL_2(\Z[i])$-equivalent to $\Phi^\iota$ give rise to curves definable
over $\Q$.

 To see this we note that by their very
construction the period matrices $Z_\Phi$ lies in a certain real
$3$-dimensional cycle in $\hh_2$ considered by  Shimura 
[Sh]. Namely, the cycle defined by 
$$
Z \in \hh_2, \qquad 
\left(\begin{array}{cc}
 0&1\\-1&0
\end{array}\right) Z=
-\overline{Z}\left(\begin{array}{cc} 0&1\\-1&0\end{array}\right) \ .
$$
If $A_Z$ is the complex abelian surface
corresponding to such a $Z$ then there is an isomorphism 
$$
\lambda: A_Z \longrightarrow \overline{A_Z}, \qquad \mbox{with }
\qquad 
\overline\lambda\circ \lambda =-\mbox{\rm id}.
$$
(Applied to $Z_\Phi$ this yields the fact that the curves $C_\Phi$
and $C_{\Phi^\iota}$ are both isomorphic and complex conjugate to each
other as mentioned above). 

It follows that if $A_Z$ has no  automorphisms other than $\pm
\mbox{\rm id}$
then it  has no model  defined over its field of moduli. It is not
hard to see that this holds for $Z_\Phi$, for every $\Phi$ in the
interior of the fundamental domain. 
\section{
Examples
}
We end with an illustration of the above discussion, giving the
outcome of algorithm  when $N=43$. The
calculations were done using PARI-GP. The routines as well as the
data for all cases is available at: \newline
\qquad \texttt{
http://www.ma.utexas.edu/users/villegas
}

The reduced forms $\Phi$ of discriminant $-43$ are $(1, 1, 11)$, $(2,
1 + 2i, 6)$ and $(3, \pm1 + 2i, 4)$.

1) For $\Phi=(2,1+2i,6)$ we obtain
$$
f_\Phi(x,1)=\frac1{2}(-\,x^6+\half(-3+567\sqrt{-43})\,x^4
+\half(3+567\sqrt{-43})\,x^2+1) 
$$
Its Igusa invariants are
$$
\begin{array}{ll}
J_2=&\phantom{ -}1728012\\
J_4=&\phantom{ -} 93313728006\\
J_6=& -186622271996\\
J_8=& -2176943579975806271997\\
J_{10}=&\phantom{ -} 2176782336000000000000
\end{array}
$$
(these were calculated using  classical algorithms for invariants of a sextic 
following  Mestre [Me]).

As in the example of the introduction $D=J_{10}=2^{-12}\disc(f_\Phi)$ factors
nicely 
$$
D=(2^2\cdot3\cdot5)^{12}.
$$

This curve descends to $\Q$; here is a  model  
$$
y^2= x^6 + 24384\,x^5 + 61311\,x^4 + 585856\,x^3 + 813483\,x^2 + 3214656\,x
+ 1472877\,.
$$

\vskip .3cm
2) For $\Phi=(3,1+2i,4)$ we obtain
\begin{eqnarray*}
f_\Phi(x,1)=\frac4{3^3}(
(14\sqrt{-43} - 160)\,x^6& + (42\sqrt{-43} + 162)\,x^5 + \\
(2247\sqrt{-43} - 159)\,x^4& + 17021\,x^3 + \\
(2247\sqrt{-43} + 159)\,x^2&
 + (-42\sqrt{-43} + 162)\,x \\
+ 14\sqrt{-43} + 160)&
\end{eqnarray*}
Its Igusa invariants are
$$
\begin{array}{ll}
J_2=&\phantom{ -}14333772\\
J_4=&\phantom{ -}7393823156166\\
J_6=&\phantom{ -}3726840435157546564\\
J_8=&-312234946681873274015037\\
J_{10}=&\phantom{ -} 7355827511386641000000000000
\end{array}
$$
and
$$
D=J_{10}=(2\cdot3\cdot5\cdot7)^{12}.
$$

As explained in \S 4, since $\Phi$ corresponds to a point in the
interior of the fundamental domain the curve $C_\Phi$ has no model
over $\Q$. Alternatively, we can see this  following Mestre [Me].
When the curve has no extra automorphisms
(i.e. its only automorphisms are the identity and the hyperelliptic
involution), the 
obstruction to the curve being definable over its field 
of moduli ($\Q$ in our case) is given by a conic in $\p^2$
$$
xMx^t=0, \qquad x=(x_1:x_2:x_3) \in \p^2,
$$
where $M$ is a $3\times 3$ symmetric matrix whose entries are
certain invariants of the sextic; more precisely, the curve is
definable over its field of moduli if and only if the conic has
a rational point there. Explicitly, we have $M=(m_{i,j})$ with 
(we have actually simplified slightly the matrix given by Mestre) 
$$
\begin{array}{ll}
m_{11}=&  3J_2^3 - 160J_4J_2 - 3600J_6\\
m_{21}=&  -J_4J_2^2 + 330J_6J_2 + 160J_4^2\\
m_{31}=&  -J_6J_2^2 -840J_6J_4 - 8000J_{10}\\
m_{22}=&  -25J_6J_2^2 - 8J_4^2J_2 -120J_6J_4 - 2000J_{10}\\
m_{32}=&  67J_6J_4 + 600J_{10}J_2 + 90J_6^2\\
m_{33}=&  -33J_6^2J_2 -100J_6J_4^2 - 800J_{10}J_4
\end{array} ,
$$
where $J_2, J_4, J_6, J_8, J_{10}$ are the Igusa invariants.

In our case we have
$$
\begin{array}{ll}
m_{11}=&-21538723388574481387776\\
m_{12}=&\phantom{ -}24856361223852137345176064256\\
m_{13}=&-23971255400369899892885589544571136\\
m_{22}=&-28732882146400381994651008552571136\\
m_{23}=&\phantom{ -}27776672840855638207256856144392139100416\\
m_{33}=&-26987491534155851141341724256178812956900004096
\end{array}
$$
We easily
verify that  this conic has rational points everywhere
locally except at the primes $43$ and $\infty$; in particular, it has
no rational points. 

We should point out that the vanishing of the determinant of $M$
precisely corresponds to the curve having extra automorphisms. As with
$D$,
this determinant factors nicely
$$
\det M=
- 2^{64} \cdot 3^{38} \cdot 5^{34} \cdot 7^{28} \cdot 19^4 \cdot 29^2 \cdot
37^2 \cdot 43\ .
$$

Finally, let $p$ be a prime which splits in $K=\Q(\sqrt{-43})$ as
$p=\pp\overline \pp$. The reduction of the curve $C_\Phi$ modulo $\pp$
gives a smooth curve $\overline C$ of genus $2$ over $\F_p$. We have
verified that for all primes in the range $167\leq p<10000$ such that
$4p=a^2+43$ for some $a\in \N$ the curve $\overline C$ or its
quadratic twist attains the maximum  number of points  possible,
namely $p+1+2\left[\sqrt{2p}\;\right]$ (an improvement on Weil bounds
due to Serre).

%
%

\end{document}